\documentclass[12pt, reqno]{amsart}

\usepackage{amssymb,latexsym,amsmath,amsfonts}
\usepackage{mathrsfs}
\usepackage{graphicx}
\usepackage[usenames]{color}
\usepackage{hyperref}
\usepackage{comment}
\usepackage{enumitem}

\definecolor{DPurple}{rgb}{0.46,0.2,0.69}

\numberwithin{equation}{section}

\allowdisplaybreaks

\theoremstyle{definition}

\theoremstyle{remark}

 \theoremstyle{plain}



\begin{document}

\title[Implications of an affirmative solution to the Lindenstrauss Problem]{Implications of an affirmative solution to the \\[1.4mm] Lindenstrauss Problem}

\parindent=0mm \vspace{.2in}
\author{M.A. Sofi}
\address{Department of Mathematics, University of Kashmir, Srinagar-190006, India.}
\email{aminsofi@gmail.com}

\keywords{Banach space, Hilbert space, Lipschitz retract, coarse map, stable. }
\subjclass[2010]{Primary: 46B03, 46B20, Secondary 46C15, 28B20.}

\begin{abstract}
The question regarding the location of Banach spaces inside their biduals has been investigated and answered reasonably satisfactorily in the linear theory of Banach spaces. Thus, for instance, whereas it is known that a dual Banach space is complemented inside its bidual, the space $c_0$ is not! However, it turns out that $c_0$ is a Lipschitz retract of $\ell_\infty$, the bidual of $c_0$. In his famous paper of 1964, Lindenstrauss asked if each Banach space is a Lipschitz retract of its bidual. In this short note, we show how to relate the Lindenstrauss problem (LP) to certain other important and well-known questions that remain open in the Lipschitz theory of Banach spaces and how these latter questions may be settled in the affirmative under the assumption of (LP) having a positive solution. 

\end{abstract}
\maketitle

\section{Introduction}\label{S:intro}

\parindent=0mm \vspace{.0in}
It is folklore, thanks to the Hahn Banach theorem, that $X$ embeds linearly and isometrically as a (closed) subspace of $X^{**}$ under the canonical embedding and that under this inclusion, Michael’s selection theorem informs us that $X$ is positioned there as a (continuous) retract of $X^{**}$. On the other hand, the principle of local reflexivity yields $X$ as a locally complemented subspace of  $X^{**}$ (See Definition 2.1). It can be proved [11] (see also [17]) that a subspace of a Banach space $X$ of which it is a Lipschitz retract is also locally complemented in $X$.

\parindent=0mm \vspace{.1in}
 In the Lipschitz geometry of Banach spaces, an old and important problem posed by J. Lindenstrauss in 1964 [14] asks whether a Banach space $X$ is always a Lipschitz retract of $X^{**}$, the bidual of $X$. As is often the case in situations involving the Lipschitz structure of Banach spaces, the question remains open in separable Banach spaces whereas a counterexample in the nonseparable case was found following an important work of Kalton [12] in 2012. In fact, Kalton’s construction is based on the non-existence of a uniform/Lipschitz right inverse of the quotient map: $\ell_{\infty}\to \ell_{\infty}/c_0$ which can also be used to construct an example of an $L_1$-predual which is not an absolute Lipschitz retract (ALR). In a subsequent development, it’s shown that such spaces can also be located within the class of (nonseparable) Schur spaces ([8]). 
 


%
%
%

\parindent=0mm \vspace{.1in}
It can be shown (see [11]) that the Lindenstrauss problem is equivalent to the following problem which is motivated by the above discussion:

\parindent=0mm \vspace{.1in}
Problem: Is it true that every locally complemented subspace of a Banach space is a Lipschitz retract.

\parindent=0mm \vspace{.1in}
We shall relate (LP) to an open problem involving the Hahn Banach extension operator constant on the one hand and to the question of stability of certain nonlinear maps in the setting of a famous theorem of Figiel on the stability of an $\varepsilon$-isometry on the other, and to point out the consequences of this relationship.

\section{Notation and background}

\parindent=0mm \vspace{.0in}
We shall throughout use the symbols $X, Y, Z$ to denote a Banach space whereas $M , N$ shall be used for metric spaces.

\parindent=0mm \vspace{.1in}
We shall also have occasion to use the following notation:

\parindent=0mm \vspace{.05in}
$\mathbb R^n$: Euclidean space of dimension $n$.

\parindent=0mm \vspace{.05in} 
$S^n$: The unit sphere in $\mathbb R^n$.

\parindent=0mm \vspace{.05in}
$B_n$: Unit ball in  $\mathbb R^n$.


\parindent=0mm \vspace{.1in}
{\bf Definition 2.1:} A subspace $Y$ of a Banach space $X$ is said to be locally complemented in $X$ if there exists $c>0$ such that for each finite dimensional subspace $M$ of $X$, there exists a continuous linear map $f:M\to Y$ with $\|f\|\le c$ and $f(x)=x$ for all $x\in M\cap Y$. In this case, we say that $Y$ is $c$-complemented in $X$. 

\parindent=0mm \vspace{.1in}
{\bf Definition 2.2:} Given a set $A$, $B\subset A$ and a mapping $f:A\to B$, we say that $f$ is a retraction onto $B$ if $f(x)=x, $ for each  x in $B$. In this case, $B$ is called a retract of $A$. Continuous, uniform, Lipschitz retracts may be defined  analogously, depending upon whether $f$ is a continuous, uniformly continuous or a Lipschitz mapping acting between topological/metric spaces.

\parindent=0mm \vspace{.1in}
{\bf Definition 2.3:} A map $f:M\to N$ acting between metric spaces $M$ and $N$ is said to be Lipschitz if for some $ k \ge 1$, we have 
\begin{align*}
\|f(x)-f(y)\|\le k\|x-y\|,\quad x,y\in M,k\ge 1.
\end{align*}

\parindent=0mm \vspace{.1in}
 For $k=1$, we say that $f$ is a nonexpansive (NE) mapping. 
 
 \parindent=0mm \vspace{.1in}
 The space of all Lipschitz maps $f : M \to \mathbb R$  shall be denoted by $Lip (M)$  whereas we shall mostly consider the (sub) space  $Lip_0 (M)$ of  $Lip (M)$ consisting of functions vanishing at a distinguished point, say $\theta \in M$. It is easily checked that under pointwise operations,$Lip (M)$ is semi-normed (and hence $Lip_0 (M)$  is a Banach) space when equipped with the (semi) norm: 
$$\|f\| = \sup_{x \ne y} \dfrac{|f(x) -f(y)}{d(x,y)}.$$

Given a metric space $(M, d),$  there is associated a Banach space $\mathcal F(M)$ the Lipschitz free Banach space of $M,$ which allows a linearization of the metric structure of $M$ by the Lipschitz properties of $M$ being transferred as linear topological properties of $\mathcal F(M)$.  The space  $\mathcal F(M)$  is defined to be the closed linear span of the set $\left\{ \delta_x : x \in M\right\}$  in $Lip_0(M)^*$ where $\delta_x(f)=f(x), f \in Lip_0(M)$  and the norm on $\mathcal F(M)$  is the one induced by $Lip_0(M)^*$    An equivalent description of the norm on $\mathcal F(M)$  is provided by looking at $\mathcal F(M)$  as the completion of the space of all measures on M of finite support with respect to the norm:
$$\mu_{\mathcal F(M)} = \inf \left\{\sum_{i=1}^{n} |\alpha_i| d(x_i ,y_i) :\mu =\sum_{i=1}^{n} \alpha_i(\delta(x_i) -\delta(y_i)) \right\}.$$

Also, it is easy  to see that the map $\delta_M : M\to \mathcal F(M)$  given by $\delta(x)=\delta_x$ embeds M in an isometric manner as a subspace of  $\mathcal F(M)$.  Sometimes, to emphasise the metric space M we shall write  $\delta_M$ for $\delta $  Further, given a Lipschitz map $f:M \to N$ there exists a bounded linear map $\tilde{f}: \mathcal F(M) \to \mathcal F(N)$  such that $\delta_N \circ f = \tilde{f}\circ \delta_M$.  It follows that if M and N are Lipschitz isomorphic, then $\mathcal F(M)$ and $\mathcal F(N)$  are linearly isomorphic. Also, in the case of M being a Banach space X, the identity map on X extends to a bounded linear map $ \beta: \mathcal F(X) \to X$  , the so called barycentric map. Thus, letting f take values in a Banach space N = X and composing  $\tilde{f}$ with $\beta$   yields a bounded linear map $\tilde{f}:\mathcal F(M) \to X$  such that $\tilde{f}\circ\delta_M =f$.  

\parindent=0mm \vspace{.1in}
For a comprehensive description of free spaces and their applications in nonlinear geometry of Banach spaces, we strongly recommend [11].

\parindent=0mm \vspace{.1in}
{\bf Example 2.4(a):} Every Banach is a continuous retract of every Banach space containing it as a subspace. In particular, and as already noted, every Banach space is a continuous retract of its bidual (See[1], Proposition 1.19).

\parindent=0mm \vspace{.1in}
{(\bf b):} (L. E. J. Brouwer): $S^n$ is not a (continuous) retract of $B^n$. Equivalently, every continuous self-map on $B^n$ admits a fixed point. 

\parindent=0mm \vspace{.1in}
On the contrary, we have the following surprising theorem in the setting of infinite dimensional spaces.

\parindent=0mm \vspace{.1in}
{\bf Example 2.5 (Benjamini, Sternfeld [2]):} In every infinite dimensional Banach space, the unit sphere is a Lipschitz retract of its unit ball. 

\parindent=0mm \vspace{.1in}

{\bf Example 2.6:} The following classes of Banach spaces which are well known to be complemented subspaces of their biduals provide easy examples of Banach spaces $X$ verifying (LP) where $X$ is a

\parindent=0mm \vspace{.1in}
a. Reflexive Banach space. More generally, 

\parindent=0mm \vspace{.1in}
b.	A dual Banach space.

\parindent=0mm \vspace{.1in}
c.	(ALR), in particular an injective Banach space.

\parindent=0mm \vspace{.1in}
The following list provides some less known and nontrivial examples of Banach spaces which are Lipschitz retracts in their bidual.

\parindent=0mm \vspace{.05in}
d.  $C(K)$ where $K$ is a compact metric space [14]

\parindent=0mm \vspace{.05in}
e. $X$, an order continuous Banach lattice [13]

\parindent=0mm \vspace{.05in}
f. $X$, separable and having an (UFDD) [13]

\parindent=0mm \vspace{.05in}

g. $K(c_0,Z)$, the space of compact linear maps taking values in an arbitrary Banach space $Z$ [13].

\parindent=0mm \vspace{.1in}
A far reaching generalisation of the example (h) above and as an important non-commutative version of the well-known fact that $c_0$ is a $2$-Lipschitz retract of $\ell_\infty$ is provided by a recent development, showing that $K(\ell_2)$ is a Lipschitz retract of its bidual: $K(\ell_2)^{**}=L(\ell_2 )$. In other words, $K(\ell_2)$) has (LP). This follows as a special case of a strengthening of this fact by a theorem of Tanaka [19] in the setting of ideals in Von Neumann algebras. The following remarkable theorem of Cheng et al.[4] provides yet another generalisation of this fact in the setting of $\ell_p$ spaces.

\parindent=0mm \vspace{.1in}
{\bf Theorem 2.7:} For $1\le p,q<\infty$, the space $K(\ell_p,\ell_q)$ of compact operators is an $8$-Lipschitz retract of $L(\ell_p,\ell_q)$.

\parindent=0mm \vspace{.1in}
{\bf Remark 2.8:} It turns out that the above conclusion cannot be strengthened to assert that the quotient map $\varphi:L(\ell_p,\ell_q)\to L(\ell_p,\ell_q)/K(\ell_p,\ell_q)$ admits a Lipschitz right inverse. Indeed, if such were the case, it would follow that the space $L(\ell_p,\ell_q)$ is uniformly homeomorphic with
\begin{align*}
L(\ell_p,\ell_q)\oplus L(\ell_p,\ell_q)/K(\ell_p,\ell_q).
\end{align*}

Further, note that $\ell_\infty$ embeds isomorphically into $L(\ell_p,\ell_q)$ whereas $c_0$ embeds into $K(\ell_p,\ell_q)$. Also, $L(\ell_p,\ell_q)$ embeds linearly and isometrically into $\ell_\infty$. Combining these observations yields that $\ell_\infty$ contains a uniformly isomorphic copy of $\ell_\infty/c_0$, contradicting Kalton’s theorem[12] that the quotient map $\ell_\infty\to \ell_\infty/c_0$  does not admit a uniformly continuous right inverse.

\parindent=0mm \vspace{.1in}
The above example motivates the question of exploring if the quotient map $X^{**}\to X^{**}/X$ admits a Lipschitz right inverse in the examples listed in 2.5.

\parindent=0mm \vspace{.1in}
As mentioned above, the main purpose of this note is to relate (LP) to certain other important open problems in the Lipschitz geometry of Banach spaces described below where certain conclusions can be shown to hold under the assumption of (LP). In the important case of separable Banach spaces, the existence of examples of Banach spaces failing the conclusion of a given theorem in the absence of (LP) as established in the subsequent sections would immediately yield counter examples to (LP) in the separable case that Kalton [13] suspects to be the case. The following well known theorem of Lindenstrauss provides a foretaste of what is to follow in this section

\parindent=0mm \vspace{.1in}
{\bf Theorem 2.9 ([14]}. See also [18], Corollary 3.19): Let a closed subspace $Y$ be a uniform retract of a Banach space $X$ such that $Y$ is a Lipschitz retract of its bidual. Then $Y$ is Lipschitz retract of $X$.

\parindent=0mm \vspace{.1in}
Thus, the existence of a closed subspace $Z$ of a Banach space $X$ which is a uniform retract but not a Lipschitz retract of $X$ would provide an example of a Banach space $Z$ failing the (LP).

\section{Results}

\parindent=0mm \vspace{.0in}

{\bf A.  Hahn Banach extension operator constant}

\parindent=0mm \vspace{.1in}
To every metric space $M$ we attach a numerical parameter $\lambda(M)$ and its vector-valued analogue $\lambda(M,Z)$ where $Z$ is a Banach space. Thus, for each pair $(M, Z)$ consisting of a metric space $M$ and a Banach space $Z$, we define the Hahn Banach extension operator constant 
\begin{align*}
\lambda(M,Z)=\sup \left\{  \lambda(A,M,Z): A\subset M\right\}
\end{align*}
where 
\begin{align*}
\lambda(A,M,Z)=\inf \left\{\|T\|:T\in Ext(A,M,Z)\right\}
\end{align*}
where $Ext(A,M,Z)$ denotes the set of (linear) extension operators: 
$$ T:Lip_0(A,Z)\to Lip_0(M,Z),\quad \text{i.e.}~~T(g)|_A =g,~~\text{for each}~g\in Lip_0(A,Z).$$

\parindent=0mm \vspace{.1in}
{\bf Definition 3.1:} Given a metric space $M$, we set $ Ext(A,M) = Ext(A,M,\mathbb R)$ and  $\lambda(M)=\lambda(M,\mathbb R)$.

\parindent=0mm \vspace{.1in}
As we shall see below, the parameter $\lambda(M)$ is closely related to the structure of the underlying metric space $M$. In the case of $M$ being a Banach space, it turns out that the finiteness of $\lambda(M)$ is equivalent to the finite dimensionality of $M$ (see [17]), whereas for a compact metrics space $M$, the finiteness of $\lambda(M)$ implies that the free space ${\mathcal F}(M)$ over $M$ has BAP [9]. 

\parindent=0mm \vspace{.1in}
In a different direction, we shall explore the nonlinear analogue of Figiel’s theorem involvingan appropriate coarse Lipschitz analogue of an $\varepsilon$-isometry acting between Banach spaces and show how under certain conditions, such maps can be approximated by a genuine Lipschitz map.

%
%
%

\parindent=0mm \vspace{.0in}


\parindent=0mm \vspace{.1in}
It is conjectured that $\lambda(M)=\lambda(M,Z)$ for each Banach space $Z$. We show that the conjecture holds for a class of Banach spaces that includes the class of  Banach spaces which are a Lipschitz retract of their biduals.

\parindent=0mm \vspace{.1in}
{\bf Theorem 3.2:} Let $M$ be a metric space and $X$ a Banach space which is a Lipschitz retract of its bidual. Then $\lambda(M)=\lambda(M,X)$.

\parindent=0mm \vspace{.1in}
{\bf Proof.}  We note that $\lambda(M)\le\lambda(M,X)$. Indeed, if $\lambda(M,X)=\infty$, there is nothing to prove. Thus, we assume that $\lambda(M,X)<\infty$. Fix $x\in Z$ such that $\|x\|=1$ and choose $g\in Z^*$ such that $g(x)=1$ with $\|g\|=1$. Let $h:\mathbb R\to X$ be the linear isometry defined by $h(t)=tx$. Let $A$ be a subset of $M$. It is easily seen that the map $E:Lip_0(A)\to Lip_0(A,X)$ given by $E(f)=hof$ defines a linear isometry. Now for a given $\epsilon>0$ there exists $T\in Ext(A,M,X)$ such that $\|T\|<\lambda (A,M,X)+\epsilon$. Then $\widehat T=goToE:Lip_0(A)\to Lip_0(M)$ defines a linear map such that $\widehat T(f)(s)=f(s)$, for each $f\in Lip_0(A)$ and $s\in A$. In other words, $\widehat T\in Ext(A,M)$ such that$\|\widehat T\|\le \|T\|< \lambda (A,M,X)+\epsilon$. This gives $\lambda (A,M)<\lambda (A,M,X)+\epsilon$. Since $A$ is arbitrary, this yields the desired inequality: $\lambda (M)\le \lambda (M,X)$.

\parindent=0mm \vspace{.1in}
For the reverse inequality, fix a subset $A$ of $M$ and assume that $\lambda (A,M)<\infty$. For each $\epsilon >0$, choose $T\in Ext(A,M)$ such that $\|T\|<\lambda (A,M)+\epsilon$. Taking conjugates gives $T^{*}:Lip_0(M)^*=\mathcal F (M)^{**}\to \mathcal F(A)^{**}=Lip_0(A)^{*}$. Since $T$ is an extension operator, it follows that $T^{*}$ is a continuous linear projection. Looking at $M$ as isometrically embedded into $\mathcal F (M)^{**}$, it follows that $T^{*}$ restricted to $M$ acts as a Lipschitz map $T^{*}\big|_M:M\to \mathcal F (A)^{**}$ such that $\left\|T^{*}|_M\right \|_{Lip(M)}\le \|T\|.$ To show that $\lambda(A,M,X)\le \lambda(A,M)$, let $f\in Lip_0(A,X)$ and consider the associated linear map $T_f:\mathcal F(A)\to X$ which extends $f$ such that $\|T_f\|\le \|f\|_{Lip_0 (A, X)}$. This gives $T_f^{**}: \mathcal F(A)^{**}\to X^{**}$ such that $T_f^{**}\big|_A=T_f\big|_A=f$. By the conditions of the theorem, there exists a Lipschitz retraction $R:X^{**}\to X$. It is readily seen that the formula: $\widehat T(f)=R \circ T_f^{**} \circ T^*\big|_M$ defines a continuous linear operator $\hat T: Lip_0 (A,X)\to Lip_0(M,X)$ with $\|\hat T(f)\|\le \|f\|_{Lip_0 (A, X)} \big(\lambda(A,M)+\epsilon  \big)$. Finally, since $\hat T(f)\big|_A=f$, it follows that $\hat T\in Ext(A,M,X)$ and that $\|\hat T\|<\lambda(A,M)+\epsilon$. Letting $\epsilon \to 0$ gives $\lambda(A,M,X)\le \lambda(A,M)$ and this yields the desired inequality.

\parindent=0mm \vspace{.1in}
{\bf Corollary 3.3:} Assume that there exists a metric space $M$ and a Banach space $Z$ for which the conjecture fails to hold, then $Z$ provides a counterexample to the Lindenstrauss problem.

\parindent=0mm \vspace{.1in}

{\bf B. 	A Lipschitz analogue of the Mazur Ulam stability}

\parindent=0mm \vspace{.1in}
In this section, the chief motivation for our study comes from an old result of Mazur and Ulam [15] asserting that a standard $(f(0)=0)$ surjective isometry $f$ acting between (real) normed spaces is a linear isometry. On the other hand, the analogous question regarding $\epsilon$-isometries was addressed in a famous theorem of Hyers and Ulam [10] that asserts that a surjective $\epsilon$- isometry (see definition below ($*$)) between Banach spaces can be approximated by a surjective (linear) isometry $g$ in the sense that the maps $f$ and $g$ stay within a finite distance from each other.

\parindent=0mm \vspace{.1in}
Regarding appropriate non-surjective analogues of these results, the first breakthrough was made by $T$. Figiel [7] who showed that in the case of the Mazur-Ulam theorem, the best one can do is that for every given standard $\epsilon$-isometry $f:X\to Y$, there exists a linear isometry $g:L(f)\to X$ of norm $1$ such that $\|gf(x)-x\|=0,~x\in X$. Here $L(f)$ denotes the closed linear span of the range of $f$. 

\parindent=0mm \vspace{.1in}
On the other hand, in recent decades certain interesting developments have been reported by a group of Chinese mathematicians involving nonsurjective $\epsilon$-isometries in the Hyers-Ulam theorem [10]. Following considerable research activity over a period of thirty years surrounding the approximation of $\epsilon$-isometries by a genuine isometry (see [3] and references therein),  Dai and his collaborators propose ([5], [6]), a Lipschitz/coarse analogue of the stability phenomenon involving the Hyers and Ulam theorem mentioned above. As it turns out, unlike the situation in the well understood case of isometry/ $\epsilon$-isometry where the stability is known to hold for certain classes of Banach spaces, the analogous situation in the coarse category as shown by Dai et al.[6] still holds, but under additional conditions including the assumption that $X$ is a Lipschitz retract of its bidual. As already pointed in Section 1, this would potentially provide new counterexamples to the Lindenstrauss problem as we shall see presently.

\parindent=0mm \vspace{.1in}
Before we discuss these issues in detail, we begin with a sample of results developed by the Chinese school on the classical case of non-surjective near isometries.

\parindent=0mm \vspace{.1in}
{\bf Theorem 3.4. ([3]): } For a given pair of Banach spaces $(X, Y)$ and $\epsilon>0$, let $f:X\to Y$ be a standard $\epsilon$-isometry:
\begin{align*}
\Big\| \|f(x)-f(y)\|-\|x-y\|  \Big\| \le \epsilon,\quad \forall~x,y\in X,~f(0)=0.\tag{*}
\end{align*}

\parindent=0mm \vspace{.0in}
Then the following statements hold:

\parindent=0mm \vspace{.1in}
(a). If $Y$ is isomorphic to a Hilbert space, there exist a bounded linear map $g:L(f)\to X$ of norm $1$ and $\alpha,\beta>0$ such that $\|g\|\le \alpha$ and $\|gf(x)-x\|\le \beta \epsilon$, for all $x\in X$.

\parindent=0mm \vspace{.1in}
 (b). Conversely, if (a) holds for all Banach spaces $X$, then $Y$ is isomorphic to a Hilbert space.

\parindent=0mm \vspace{.1in}
(c). For $X$, an injective Banach space, there exist a bounded linear map $g:L(f)\to X$ of norm $1$ and $\alpha, \beta>0$ such that $\|g\|\le \alpha$ and $\|gf(x)-x\|\le \beta \epsilon$, for all $x\in X$.

\parindent=0mm \vspace{.1in}
(d).	Conversely, if (c) holds for all Banach spaces $Y$, the $X$ is cardinality injective.

\parindent=0mm \vspace{.1in}
Here, {\it cardinality injectivite} refers to a Banach space $X$ which is complemented in each larger Banach space having the same cardinality as $X$.

\parindent=0mm \vspace{.1in}
We shall use the term universally right stable (URS) for those Banach spaces $Y$ for which the condition in Theorem 3.4(a) holds for all Banach spaces $X$ whereas $X$ shall be called a universally left stable (ULS) space if the same holds for all Banach spaces $Y$.

\parindent=0mm \vspace{.1in}
In the nonlinear category, the idea is to explore if an appropriate analogue of the above ideas involving quantitative versions of Figiel’s theorem may be realized in the Lipschitz category of Banach spaces; in other words, whether a (standard) coarse Lipschitz embedding may be approximated (in the sense of Theorem 3.4(a)) by a Lipschitz map, the sort of mappings which are used as morphisms in the Lipschitz category. 

\parindent=0mm \vspace{.1in}
{\bf Definition 3.5(a):} A mapping $f:(M,d_1)\to (N,d_2)$ acting between metric
spaces $M$ and $N$ is said to be coarse continuous  (or bornologous) if $\forall~ \epsilon>0$,  $\exists \, \delta>0$ such that
\begin{align*}
d_1(x,y)<\epsilon \Rightarrow d_2(f(x),f(y))< \delta,~\forall~x,y\in M.
\end{align*}

\parindent=0mm \vspace{.0in}
(b). $f$ is said to be a {\it coarse embedding} if both $f$ and its inverse (from the range of $f$) are coarse maps. 

\parindent=0mm \vspace{.1in}
(c). A map $f:M\to N$ is said to be {\it Lipschitz for large distances} if there exist $K,L>0$ such that $d_2(f(x),f(y) )\le Kd_1(x,y)+L,~\forall ~x,y\in M$. In particular, there exist $K,L>0$ such that $(d_2 f(x),f(y) )\le Kd_1 (x,y)$, for all $x,y\in M$ with $d_1 (x,y)\ge L$. 

\parindent=0mm \vspace{.1in}
Further, $f$ is said to be a quasi-isometric embedding(also called coarse Lipschitz embedding in the sense of [6]) if there exist $K, L, R, S > 0$ such that
$$ Rd_1 (x,y)-S\le d_2 (f(x),f(y) )\le Kd_1 (x,y)+L,~\forall~ x,y\in M.$$

\parindent=0mm \vspace{.1in}
{\bf Definition 3.6(a):} A pair $(X, Y)$ of Banach spaces is said to be {\it coarsely stable} if the following holds:

\parindent=0mm \vspace{.1in}
For every standard quasi-isometric embedding $f:X\to Y$, there exist $\alpha, \beta>0$and a Lipschitz map  $g:L(f)\to X$ with $\|g\|_{Lip} \le \alpha $ such that $\|gf(x)-x\|\le \beta,$ for all $x\in X$.

\parindent=0mm \vspace{.1in}
(b). A Banach space $X$ (resp.$Y$) is said to be universally left (resp. right) coarsely stable (ULCS)(resp. URCS) if the pair $(X,Y)$ is coarsely stable for each Banach space $Y$ (resp. $X$).

\parindent=0mm \vspace{.1in}
Unlike in the linear case where the equivalence $((a)\Leftrightarrow (b))$ in Theorem 3.4 always holds as noted above, in the nonlinear case involving equivalent conditions for coarse Lipschitz embeddings, Dai et al [6] show that the reverse implication holds under certain additional conditions including the reflexivity of $X$. The upshot of Dai’s approach to Lindenstrauss’s problem is summed up in the following interesting theorem. Here, we shall use the  adjective Hilbertian to denote a Banach space which is linearly isomorphic to a Hilbert space.

\parindent=0mm \vspace{.1in}
{\bf Theorem 3.7(a)([6]):} $X$, Hilbert space $\Rightarrow X$ is (URCS).

\parindent=0mm \vspace{.1in}
(b). Converse: $X$ reflexive and has (URCS) $\Rightarrow X$ is Hilbertian.

(c). More generally, if Lindenstrauss Problem has a positive solution, then for a Banach space $X$ with (RNP), we have: (URCS) $\Rightarrow X$ is Hilbertian. Equivalently, under these conditions, we have the equivalence:  $X$ is Hilbertian $\Leftrightarrow X$ is (URCS).

\parindent=0mm \vspace{.1in}
The main result of this section that follows next offers a refinement of Dai’s theorem given above by showing that the condition in Theorem 3.7(c) can be considerably relaxed in that the assumption of (RNP) in (c) can be dispensed with.

\parindent=0mm \vspace{.1in}
{\bf Theorem 3.8:} For a Banach space $X$, consider the following statements: isomorphic to a Hilbert space.

\parindent=0mm \vspace{.1in}
(a). $X$ is Hilbertian

\parindent=0mm \vspace{.05in}
(b). $X$ is (URS)

\parindent=0mm \vspace{.05in}
(c). $X$ is (URCS).

\parindent=0mm \vspace{.1in}
Then $(a)\Leftrightarrow (b) \Rightarrow (c)$. Further, if (LP) has a positive solution, then $(c)\Rightarrow (a)$ and, therefore, all the conditions are equivalent. 

\parindent=0mm \vspace{.1in}
{\bf Proof:} The equivalence $(a)\Leftrightarrow(b)$ follows from Theorem 3.4 whereas $(a)\Rightarrow(c)$ is showed in [5, Theorem 2]. To show that $(c)\Rightarrow(a)$, assume that $X$ is (URCS). To this end, we follow an idea of Qian [16].  Let $Y$ be a closed subspace of $X$. Choose a standard 1-1 mapping $h:Y\to B_X$ and define $f:Y\to X$ by $f(y)=y+\frac{\epsilon}{2} h(y)$. Being an $\epsilon$-isometry, $f$ is clearly a coarse Lipschitz embedding such that $L(f)=X$. By definition, since $(Y, X)$ is (URCS), there exist $\alpha>0,\beta>0$ and a Lipschitz map $g:X\to Y$ with $\|g\|_Lip\le \alpha $ such that $\|gf(y)-y\|\le \beta,$ for all $y\in Y$. We note that since for each $y\in Y, \lim_{n\to \infty} \frac{f(ny)}{n}=y$, it follows that $\lim_{n\to \infty} \frac{g(ny)}{n}=y$ for each $y\in Y$. However, an application of Banach-Alaoglu theorem yields for each $x\in X$ that the sequence $\frac{g(nx)}{n}$ has a weak$^*$ limit, say $P(x)$  in $Y^{**}$ over any free ultrafilter $\mathcal U$ on $\mathbb N$. In other words, for each  $x\in X:$
$$P(x)=w^{*}-\lim_{\mathcal U}\frac{g(nx)}{n}.$$
Note that since the dual norm on$Y^{**}$ is $\mathcal U- w^*$ lower semi-continuous, it follows that for all $x,y\in X$, we have
\begin{align*}
\big\|  P(x)-p(y)\big\|&=\left\|w^{*}-\lim_{\mathcal U}\frac{g(nx)}{n}-w^{*}-\lim_{\mathcal U}\frac{g(ny)}{n}  \right\|\\
&\le \lim_{\mathcal U}\left\|\dfrac{g(nx)-g(ny)}{n}  \right\|= \lim_{n}\left\|\dfrac{g(nx)-g(ny)}{n}  \right\|\le \alpha\|x-y\|.
\end{align*}

\parindent=0mm \vspace{.0in}
It follows that $P:X\to Y^{**}$ is a Lipschitz mapping such that $P(y)=y$ for each $y\in Y$, with $\|P\|_{Lip}\le \alpha$. Composing $P$ with the assumed retraction $r:Y^{**}\to Y$ gives that $Y$ is a Lipschitz retract of $X$. By [17], Corollary 13, it follows that $X$ is (isomorphic to) a Hilbert space.

\parindent=0mm \vspace{.1in}
Regarding the analogous statement for the dual property involving (ULCS), we have the following.

\parindent=0mm \vspace{.1in}
{\bf Theorem 3.9 (a) ([6]):} $X$ is an absolute cardinality Lipschitz retract $\Rightarrow X$ is (ULCS).

\parindent=0mm \vspace{.1in}
(b) If $X$ is a Lipschitz retract of $X^{**}$, then $X$ is (ULCS) $\Rightarrow X$ is absolute cardinality Lipschitz retract.

\parindent=0mm \vspace{.1in}
{\bf Remark: 3.10:}  It would be highly desirable if, as in the case of its linear
counterpart, the conclusions in Theorem 3.8(c) $\Rightarrow $ (a) and Theorem 3.9(b) $\Rightarrow$ (a) involving the Lipschitz/coarse category, would hold without assuming the truth of (LP). On the other hand, were it to be the case that the said implications do not hold in the absence of (LP), more so in the separable case, then we would have a much awaited counterexample to (LP), as summed up in the Corollary below.

\parindent=0mm \vspace{.1in}
{\bf Corollary 3.11:} Lindenstrauss Problem has a negative answer for Banach spaces $X$ in each of the following cases: 

\parindent=0mm \vspace{.1in}
(a). There exists a metric space $M$ and a Banach space $X$ such that the equality $\lambda(M)=\lambda(M,X)$ does not hold. 
	
\parindent=0mm \vspace{.05in}	
(b). There exists a Banach space $X$ which is universally right coarse stable (URCS) but is not Hilbertian.

\parindent=0mm \vspace{.05in}	
(c). There exists a Banach space $X$ which is universally left coarse stable (ULCS) but is not an absolute cardinality Lipschitz retract.

\parindent=0mm \vspace{.3in}

 {\bf{$*$Availability of data and material}}: Not applicable.

\parindent=0mm \vspace{.1in}
  {\bf{$*$Competing interests:}} There is no conflict of interest.

\parindent=0mm \vspace{.1in}
{\bf{$*$ Acknowledgement:}} The author would like to thank the referee for his/her suggestions that have helped improve the presentation of the material.

\end{document}